\documentclass{IEEEtran}
\usepackage{graphicx}
\usepackage{subfigure}

\usepackage[colorlinks]{hyperref}

\usepackage{amsmath}
\usepackage{amssymb}
\usepackage{cite}

\usepackage[ruled]{algorithm2e}
\usepackage{algorithmic}
\usepackage{multirow}

\usepackage{footnote}

\usepackage{array}
\makeatletter
\newcommand{\thickhline}{%
    \noalign {\ifnum 0=`}\fi \hrule height 1pt
    \futurelet \reserved@a \@xhline
}
\newcolumntype{"}{@{\hskip\tabcolsep\vrule width 1pt\hskip\tabcolsep}}
\makeatother

\begin{document}
%
\title{Port Selection for Fluid Antenna Systems\thanks{The work is supported in part by EPSRC under grant EP/T015985/1.}}
%
%
%

\author{Zhi Chai, 
            Kai-Kit Wong,~\IEEEmembership{Fellow,~IEEE,} 
            Kin-Fai Tong,~\IEEEmembership{Senior Member,~IEEE,}\\
            Yu Chen,~\IEEEmembership{Member,~IEEE,} 
            and Yangyang Zhang

\thanks{Z. Chai, K. K. Wong and K. F. Tong are with the Department of Electronic and Electrical Engineering, University College London, London WC1E 7JE, United Kingdom. Corresponding author: $\rm kai\text{-}kit.wong@ucl.ac.uk$.}
\thanks{Y. Chen is with Beijing University of Posts and Telecommunications, China.}
\thanks{Y. Zhang is with Kuang-Chi Science Limited, Hong Kong SAR, China.}
}

%
%

\markboth{Submitted to IEEE Communications Letters, 2022}%
{6G Wireless System}
%



\maketitle

\begin{abstract}
Fluid antenna system promises to obtain enormous diversity in the small space of a mobile device by switching the position of the radiating element to the most desirable position from a large number of prescribed locations of the given space. Previous researches have revealed the promising performance of fluid antenna systems if the position with the maximum received signal-to-noise ratio (SNR) is chosen. However, selecting the best position, referred to as {\em port selection}, requires a huge number of SNR observations from the ports and may prove to be infeasible. This letter tackles this problem by devising a number of fast port selection algorithms utilizing a combination of machine learning methods and analytical approximation when the system observes only a few ports. Simulation results illustrate that with only $10\%$ of the ports observed, more than an order of magnitude reduction in the outage probability can be achieved. Even in the extreme cases where only one port is observed, considerable performance improvements are possible using the proposed algorithms.
\end{abstract}


%
\IEEEpeerreviewmaketitle

\section{Introduction}
\IEEEPARstart{D}{iversity} and multiplexing gain are interconvertible by advanced coding and signal processing schemes \cite{Tse-2003}. This has provided the foundation for multiple-input multiple-output (MIMO) to flourish in recent-generations mobile communications. The introduction of massive MIMO means that in 5G, there are $64$ antennas at a base station (BS) \cite{Larsson-2019,Huawei-2019}. However, the same rise in the number of antennas at a user equipment (UE) is not foreseen for the reason that the space of the UE is limited even though the antenna size is getting smaller when moving up the frequency bands. Furthermore, the challenge of incorporating multiple antennas is that the common practice is to ensure antenna spacing of at least half a wavelength for sufficient diversity and minimal mutual coupling.

Recently, it has emerged that space diversity can be obtained by a very flexible, fluidic conductive antenna structure which is referred to as {\em fluid antenna} \cite{Wong-2020ell}. Such flexible antennas should not come as a surprise following the successes by Mitsubishi Electric showing a radiation efficiency of $70\%$ using seawater for antennas \cite{Mitsubishi} and the work by Xing {\em et al.} achieving $360$-degree beam-steering using a saltwater-based antenna \cite{Xing-19}. In addition, reconfigurable fluid antennas have been designed, resulting in a wide range of functionalities, e.g., \cite{Ohta-2013,Dey-2016,Tong-2017,Singh-2019,Shen-2020apmc}. A contemporary survey on this topic can be found in \cite{Huang-2021access}. 

{\em So, what does it mean to have fluid antenna at a UE?}

One outcome is that we can now have a position-switchable antenna which can maximize its received signal-to-noise ratio (SNR) by choosing the best position (i.e., port selection) in a given space. The beauty is that diversity comes essentially by a single antenna in the manner of selection combining without the concern of mutual coupling. In \cite{Wong-2020twc}, such position-tuneable fluid antenna system was analyzed and it was revealed that a fluid antenna system could match the outage performance of a maximal-ratio combining (MRC) receiver with many spatially uncorrelated antennas if the number of switchable positions, referred to as `ports', was sufficiently large even if the space is small at the UE. The capacity gain using fluid antenna for a point-to-point system was also recently analyzed in \cite{Wong-2020cl}.

The working principle of fluid antenna is simple. It skims through a large collection of fading envelopes from the ports and switches the radiating element to the one with the highest peak for the maximum SNR. To have a large gain, the number of ports, $N$, should be large, which can be practically realized using the surface-wave based architecture in \cite{Shen-2020apmc}. Estimating the SNRs for the very large number of ports, unfortunately, is infeasible, if not impossible, as it incurs delay in switching to each port for SNR observation.\footnote{Switching materials from one place to another causes delay. With the use of nano-pumps and higher frequency bands, the diameter of the micro-fluidic system for the fluid antenna would be less than $1{\rm mm}$, and the response time would be in the sub-millisecond range \cite{Convery-2019} and such delay can be negligible. However, port switching many times would cause unbearable delays.} This motivates our work to tackle the port selection problem in which only a very small number of ports of the fluid antenna are observed. 

In particular, the aim of this letter is to develop port selection algorithms that can approach the performance of optimal port selection when only the SNRs of a few ports are observed. This is considered possible because of the strong spatial correlation among the ports in a tight space. We will present a number of machine learning based algorithms for port selection of fluid antenna systems. The results are indispensable in realizing the benefits of fluid antenna in practice. {\color{blue}Machine learning methods have recently been applied in antenna selection problems \cite{Joung-2016,He-2018} but our problem has two unique features, namely (1) an extremely large number of ports for selection based on limited observations, and (2) strong channel correlation due to high density of the ports, which have not been considered before.\footnote{\color{blue}The features suggest that the channel envelope follow certain patterns as a result of the correlation. Machine learning, being a well-known tool to find subtle patterns, becomes a natural choice to our problem.} 

Our contributions can be summarized as follows:
\begin{itemize}
\item Using a mixture of machine learning methods including the new framework, Smart, `Predict and Optimize' (SPO) \cite{Grigas-2021}, we devise algorithms that infer the best port based on only a few port observations. Results demonstrate that great reduction in outage probability can be obtained even with {\em only one port observation} using the methods.

\item Additionally, fluid antenna with learning based port selection outperforms considerably multiple antenna systems with best antenna selection occupying the same space. This provides a strong justification to the practicality of fluid antenna systems as a new means for spatial diversity, contrasting to systems with multiple fixed antennas.
\end{itemize}}


\section{Fluid Antenna System}\label{sec:fas}
\subsection{System Model}\label{ssec:model}
We consider a point-to-point system in which a transmitter using a standard antenna is sending information to a mobile receiver equipped with a fluid antenna. The fluid antenna can switch the location of its radiating element to one of the $N$ preset locations (or ports) evenly distributed along a linear space of length, $W\lambda$ in which $\lambda$ is the wavelength. The delay for port switching is ignored. Under flat fading, the received signal at the $k$-th port (time index omitted) is given by
\begin{equation}\label{eqn:yk}
z_k=g_k s+\eta_k,
\end{equation}
where $g_k$ denotes the complex channel at the $k$-th port, which is complex Gaussian distributed with zero mean and variance of $\sigma^2$, $\eta_k$ is the zero-mean complex Gaussian noise at the $k$-th port with variance of $\sigma_\eta^2$, and $s$ represents the information symbol. The received average SNR at each port is found as
\begin{equation}\label{eqn:snr}
\Gamma=\sigma^2\frac{{\rm E}[|s|^2]}{\sigma^2_\eta}\equiv\sigma^2\Theta,~\mbox{where }\Theta\triangleq\frac{{\rm E}[|s|^2]}{\sigma^2_\eta}.
\end{equation}

The channels {\color{blue}$\{g_k\}_{\forall k}$} are considered to be spatially correlated because they can be arbitrarily close to each other. For small $W$, it is expected that {\color{blue}$\{g_k\}_{\forall k}$} are highly correlated. To characterize the correlation, we measure the displacement of the $k$-th port from the first port as
\begin{equation}\label{eqn:dn}
d_k=\left(\frac{k-1}{N-1}\right)W\lambda,~\mbox{for }k=1,2,\dots,N.
\end{equation}
In this letter, the amplitude of the channel at each port, $|g_k|$, is Rayleigh distributed, with ${\rm E}[|g_k|^2]=\sigma^2$. With 2-D isotropic scattering and assuming isotropic ports, the cross-correlation function of the channel ports satisfies \cite{Stuber-2002}
\begin{equation}\label{eqn:space-correlation}
\phi_{g_kg_\ell}(d_k-d_{\ell})
=\frac{\sigma^2}{2}J_0\left(\frac{2\pi(k-\ell)}{N-1}W\right),
\end{equation}
where $J_0(\cdot)$ is the zero-order Bessel function of the first kind. 

\subsection{Port Selection}\label{ssec:port}
\begin{figure}[]
\begin{center}
\includegraphics[width=\linewidth]{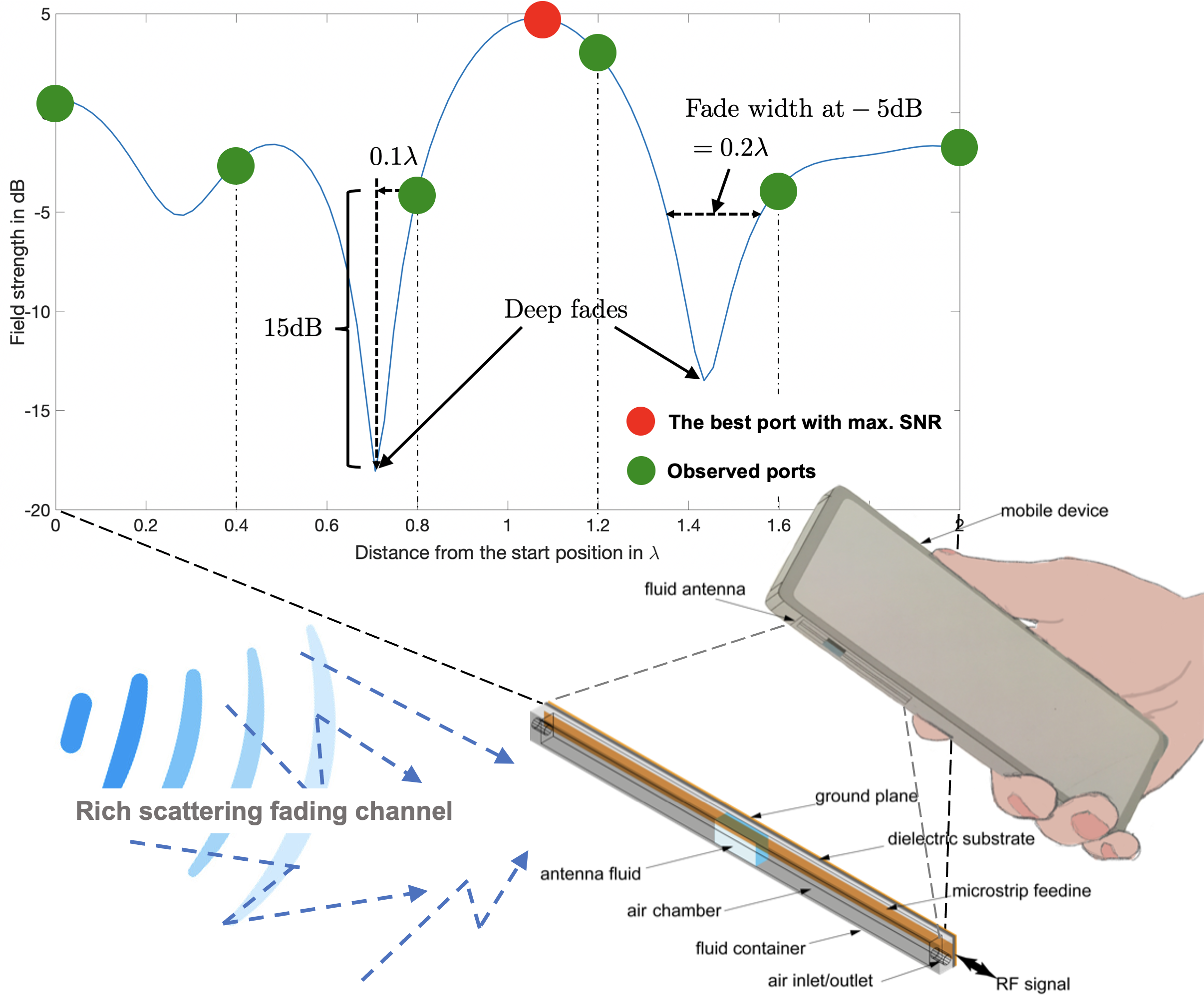}
\caption{Port selection of a fluid antenna knowing only a subset of ports.}\label{fig:portselect}
\end{center}
\end{figure}

To optimize the system performance, the aim is to switch the radiating element to the best port for the maximum received SNR. Mathematically, this is expressed as
\begin{equation}\label{eq:portselect-1}
k_{\rm opt}=\arg_k\max\left\{|g_1|,|g_2|,\dots,|g_N|\right\},
\end{equation}
where $k_{\rm opt}$ gives the index of the optimal port. The search is straightforward if the receiver knows all $\{|g_k|\}$. Nevertheless, in practice, the number of ports, $N$, for a fluid antenna can be very large\footnote{Note that if the surface-wave based architecture in \cite{Shen-2020apmc} is used, the number of ports for the fluid antenna will only be limited by the resolution of the digital control of the pump. Hence, $N$ can be extremely large.} and estimating $|g_k|$ for all the ports is infeasible. Therefore, in this letter, we consider a more practical scenario where only a small subset of the ports are observed. As such, the problem (\ref{eq:portselect-1}) becomes (see Fig.~\ref{fig:portselect})
\begin{equation}\label{eq:portselect-2}
{\color{blue}k^*=\arg_k\max\left\{\{|g_k|\}_{k\in{\cal K}},\{|\tilde{g}_k|\}_{k\in{\cal U}}\right\},}
\end{equation}
where ${\cal K}$ denotes the set of the indices of known channel gains by direct observation, ${\cal U}$ is the set of the indices of estimated channel gains by other means, and $\tilde{g}_k$ represents an estimated channel at the $k$-th port if it has been estimated.

The problem (\ref{eq:portselect-2}) suggests that the unobserved ports should be estimated somehow so that port selection can be conducted. In the sequel, we develop algorithms for estimating the best port directly or the unobserved ports for port selection.

\section{Proposed Algorithms}\label{sec:alg}
In this section, we present several algorithms that solve (\ref{eq:portselect-2}) suboptimally. All of them exploit the spatial correlation over the ports of the fluid antenna, assuming a linear structure.

\subsection{Analytical Approximation (AA)}
Given the correlation structure (\ref{eqn:space-correlation}), it makes sense to model the channels at the $N$ antenna ports by 
\begin{equation}\label{eqn:g-model}
\left\{\begin{aligned}
g_1&=\sigma x_0 +j\sigma y_0\\
g_k&=\sigma\left(\sqrt{1-\mu_k^2}x_k+\mu_N x_0\right)\\
&\quad\quad\quad\quad+j\sigma\left(\sqrt{1-\mu_k^2}y_k+\mu_k y_0\right), k=2,\dots,N,
\end{aligned}\right.
\end{equation}
in which $x_0,\dots,x_N,y_0,\dots,y_N$ are all independent Gaussian random variables with zero mean and variance of $1/2$, and $\{\mu_k\}$ are the cross-correlation parameters {\color{blue}given by \cite{Stuber-2002}}
\begin{equation}\label{eqn:mu-condition}
\mu_k=J_0\left(\frac{2\pi(k-1)}{N-1}W\right),~\mbox{for }k=2,\dots,N.
\end{equation}
The objective here is to construct an estimate of $|g_\ell|$ for all $\ell\in{\cal U}$ given the known (or observed) channels $\{|g_k|\}_{k\in{\cal K}}$. We understand that the probability density function (pdf) of $|g_\ell|$ for $\ell\in{\cal U}\backslash\{1\}$ conditioned on the known channels is 
\begin{align}
p_{|g_\ell|}&(r|\{|g_k|\}_{k\in{\cal K}})\notag\\
&\stackrel{(a)}{=}p_{|g_\ell|}(r|x_0,y_0)\notag\\
&\stackrel{(b)}{=}\frac{2r^2}{\sigma^2(1-\mu^2_\ell)}e^{-\frac{r^2+\sigma^2\mu^2_\ell(x_0^2+y_0^2)}{\sigma^2(1-\mu^2_\ell)}}
I_0\left(\frac{2\mu_\ell r\sqrt{x_0^2+y_0^2}}{\sigma(1-\mu^2_\ell)}\right),\label{eq:cpdf}
\end{align}
where $(a)$ uses the fact that the correlation between the ports is connected only by the random variables, $x_0$ and $y_0$ and the respective cross-correlation parameter, $(b)$ comes from the fact that the conditional pdf is Rician distribution with the mean determined by $x_0$ and $y_0$, and $I_0(\cdot)$ is the zero-order modified Bessel function of the first kind. As a result, it is possible to provide an estimate for $|g_\ell|$ by generating a random sample based on the conditional pdf in (\ref{eq:cpdf}). Using this technique, the AA method can estimate all the unobserved ports easily.

\subsection{Long Short-Term Memory (LSTM)}
In order to learn the correlation structure between the ports and exploit it for the estimation of $\{|g_\ell|\}_{\ell\in{\cal U}}$, machine learning techniques are useful and in particular, LSTM is a special type of recurrent neural network that memorizes information in a time-series manner. LSTM has been proved to be effective in many communications systems such as wireless caching \cite{Tao-2021}. Here, we apply LSTM to treat the channel gains over the ports in space as a time series (a space series to be precise in our case). In this way, the forecasting capability of LSTM can be used to estimate the unobserved ports. Supervised learning is adopted to train the LSTM model with the mean-square-error (MSE) as the {\color{blue}training stage} loss function, i.e.,
\begin{equation}\label{eq:loss-mse}
{\cal L}_{\rm mse}={\rm E}\left[\frac{1}{|{\cal U}|}\sum_{\ell\in{\cal U}}(|g_\ell|-|\tilde{g}_\ell|)^2\right],
\end{equation}
where $\tilde{g}_\ell$ is the estimate as the output of the LSTM model, $g_\ell$ is the true channel gain known in the dataset for training, and the expectation is taken over all the examples in the dataset. The detailed parameters of the LSTM model will be provided in Section \ref{sec_results} when simulation results are discussed.

\subsection{SPO}
While estimating the unobserved ports before port selection is a logical thing to do, this does not always return the best estimate of the best port. A better way is to base the estimation on the decision one wants to make. This was tackled in \cite{Grigas-2021} where the SPO framework which measured the decision error induced by a prediction was proposed. In a nutshell, SPO is designed to perform estimation (or prediction) with emphasis on the outcome of the final optimization. Specifically, for our port selection problem, this leads to the loss function
\begin{equation}\label{eq:loss-spo}
{\cal L}_{\rm spo}(\tilde{\boldsymbol{g}},\boldsymbol{g})=\boldsymbol{g}^T\omega^*(\tilde{\boldsymbol{g}})-z^*(\boldsymbol{g}),
\end{equation}
where $\boldsymbol{g}\triangleq [|g_1|\cdots |g_N|]^T$, $\tilde{\boldsymbol{g}}\triangleq [|\tilde{g}_1|\cdots |\tilde{g}_N|]^T$, the superscript $T$ denotes the transpose of a vector, 
\begin{equation}\label{eq:bport}
\omega^*(\boldsymbol{g})=\arg\max_{\boldsymbol{x}\in{\cal X}}\boldsymbol{g}^T\boldsymbol{x},
\end{equation}
where ${\cal X}$ is the set of all possible $N$-dimensional vectors that have zeros in all dimensions except one being unity, and 
\begin{equation}
z^*(\boldsymbol{g})=\max_{\boldsymbol{x}\in{\cal X}}\boldsymbol{g}^T\boldsymbol{x}.
\end{equation}
The second term of (\ref{eq:loss-spo}) gives the channel gain of the best port while the first term finds the channel gain of the estimated port. Intriguingly, in (\ref{eq:loss-spo}), there is no restriction on the accuracy of the estimated unobserved ports. In other words, as long as $\omega^*(\tilde{\boldsymbol{g}})$ outputs the correct port, it is as good as perfect even if some estimates deviate greatly from the true values.

The minimization (\ref{eq:loss-spo}) is however NP-hard and the following approximation of ${\cal L}_{\rm spo}(\tilde{\boldsymbol{g}},\boldsymbol{g})$ has been proposed \cite{Grigas-2021}:
\begin{multline}\label{eq:loss-spo+}
{\cal L}_{{\rm spo}+}(\tilde{\boldsymbol{g}},\boldsymbol{g})=(\boldsymbol{g}^T-2\tilde{\boldsymbol{g}}^T)\omega^*(\boldsymbol{g}^T-2\tilde{\boldsymbol{g}}^T)\\
+2\tilde{\boldsymbol{g}}^T\omega^*(\boldsymbol{g})-\boldsymbol{g}^T\omega^*(\boldsymbol{g}),
\end{multline}
in which $\tilde{\boldsymbol{g}}=f(\boldsymbol{a})$ is obtained by some mapping function $f$ over the feature vector $\boldsymbol{a}$ which basically contains the channel gains of all the observed ports. In this letter, we adopt a linear mapping so that $\tilde{\boldsymbol{g}}={\bf B}\boldsymbol{a}$ where ${\bf B}$ is an $N\times|{\cal K}|$ matrix. Note that in so doing, (\ref{eq:loss-spo+}) combines the process of estimation and optimization together. Despite this, the minimization of (\ref{eq:loss-spo+}) is still not solvable because $\omega^*(\cdot)$ involves maximization over a discrete set. Fortunately, according to \cite{Jaggi-2011}, we can replace the feasible set ${\cal X}$ by its convex hull without changing the optimal value of (\ref{eq:bport}). Hence, we can write
\begin{equation}\label{eq:bport-2}
\omega^*(\boldsymbol{g})=\arg\max_{\boldsymbol{x}\in{\rm conv}({\cal X})}\boldsymbol{g}^T\boldsymbol{x}.
\end{equation}
With (\ref{eq:bport-2}), (\ref{eq:loss-spo+}) can be minimized using a subgradient method, with the subgradient of (\ref{eq:loss-spo+}) derived as
\begin{equation}\label{eq:subgrad}
s(\tilde{\boldsymbol{g}},\boldsymbol{g})=2\left[\omega^*(\boldsymbol{g})-\omega^*(2\tilde{\boldsymbol{g}}-\boldsymbol{g})\right].
\end{equation}

Overall, what we need to obtain is the mapping ${\bf B}$, which eventually gives out the port selection solution $\omega^*({\bf B}\boldsymbol{a})$. This can be done by solving
\begin{equation}
\min_{{\bf B}\in\mathbb{R}^{N\times|{\cal K}|}}\sum_{i=1}^Q{\cal L}_{{\rm spo}+}({\bf B}\boldsymbol{a}^{(i)},\boldsymbol{g}^{(i)}),
\end{equation}
where $Q$ is the batch size or the number of training examples in the dataset, $\boldsymbol{a}^{(i)}$ denotes the feature vector of the $i$-th sample, and $\boldsymbol{g}^{(i)}$ is $i$-th labelled sample. A stochastic gradient descent algorithm is proposed for updating ${\bf B}$ (see Algorithm \ref{alg:2}).

\begin{algorithm}[]
\caption{Stochastic Gradient Descent (SGD)}\label{alg:2}
\begin{algorithmic}[1]
\STATE {\bf Initialize} ${\bf B}_1\in\mathbb{R}^{N\times|{\cal K}|}$ randomly
\STATE {\bf Set} $i=t=1$
\STATE {\bf Repeat}
\STATE \quad\quad{\bf While} {$i\le Q$} {\bf compute}
\begin{align*}
\tilde{\omega}_t^{(i)}&=\omega^*(2{\bf B}_t\boldsymbol{a}^{(i)}-\boldsymbol{g}^{(i)})\\
\tilde{\bf G}_t^{(i)}&=(\omega^*(\boldsymbol{g}^{(i)})-\tilde{\omega}_t^{(i)})(\boldsymbol{a}^{(i)})^T\\
i&=i+1
\end{align*}
\STATE \quad\quad{\bf end While}
\begin{align*}
{\bf G}_t&=\frac{1}{Q}\sum_{i=1}^Q\tilde{\bf G}_t^{(i)}\quad\quad\quad\quad\quad\quad\quad\quad\quad\quad\quad\quad\quad\quad\quad\\
{\bf B}_{t+1}&={\bf B}_t-\alpha{\bf G}_t,~\mbox{for some }\alpha>0\quad\quad\quad\quad\quad\quad\quad\quad\quad\quad\quad\quad\quad\quad\quad\\
t&=t+1\quad\quad\quad\quad\quad\quad\quad\quad\quad\quad\quad\quad\quad\quad\quad
\end{align*}
\STATE {\bf Until convergence}
\end{algorithmic}
\end{algorithm}

\subsection{Other Variants}
It is possible to mix the above approaches to perform better. Given a dataset of $Q$ labelled examples, we can split it into two datasets, half for training and another half for testing SPO. The output for testing of SPO produces a feature set that can be used to train and test LSTM if SPO and LSTM work in a concatenated manner. Two-third of the feature set can be used to train LSTM while the rest is used to test it. Another option is to use AA to preprocess the labelled dataset before it can be used for SPO. We refer to these approaches as SPO+LSTM and AA+SPO+LSTM, respectively. {\color{blue}Table \ref{table:com} provides the computational complexity of the different schemes using the big-O notations, where we focus on how the complexity scales with the parameters of the fluid antenna, $|{\cal K}|$ and $|{\cal U}|$.}

\begin{table}[!h]
\centering
{\color{blue}
\caption{Online computational complexity for the schemes} 
\begin{tabular}{c|l} 
\thickhline 
Method & Details\\ 
\hline\hline 
SPO & $O(mn)$ \\\hline 
LSTM & $O(n)$ \\\hline
AA & $O(n)$ \\\hline
SPO+LSTM & $O(mn+n)$ \\\hline
AA+SPO+LSTM & $O(mn+n)$ \\ [0.1ex] 
\thickhline 
\multicolumn{2}{l}{\footnotesize $O(\cdot)$ denotes the big-O notation.} \\
\multicolumn{2}{l}{\footnotesize $m=|{\cal U}|$ represents the number of unobserved ports.} \\
\multicolumn{2}{l}{\footnotesize $n=|{\cal K}|$ represents the number of observed ports.} \\
\end{tabular}\label{table:com}}
\end{table}


\section{Simulation Results}\label{sec_results}
In this section, we provide simulation results to assess the proposed algorithms. Rayleigh fading was assumed in all the simulations of the channel envelopes and a linear structure of the fluid antenna was considered. {\color{blue}Also, the observed ports in ${\cal K}$ were chosen so that they spread evenly over the space of $W\lambda$, as this was the most effective way to exploit the correlation. For example, if $|{\cal K}|=1$, then the $25^{\rm th}$ port will be chosen as the observed one. If $|{\cal K}|=5$, then the $1^{\rm st}$, $12^{\rm th}$, $25^{\rm th}$, $38^{\rm th}$ and $50^{\rm th}$ ports will be the observed ones, and so on.\footnote{\color{blue}We assume that the ports are numbered from $1$ to $N$ from one end of the fluid antenna to another end.} With a given number of observed ports, we use `Reference' to represent the system choosing the best port out from the observed ports.\footnote{\color{blue}Notice that `Reference' is not the same as the system selecting the best antenna with the number of fixed uncorrelated antennas equalling the number of observed ports. The reason is that for `Reference', the observed ports are usually correlated. Hence, its performance is expected to be inferior to the best antenna selection system with uncorrelated antennas.} Another benchmark is the best antenna selection system with the maximum number of uncorrelated antennas allowed in the given space. If $W=0.5$, then there can be two fixed antennas with $\frac{\lambda}{2}$ apart, while for $W=2$, $5$ antennas can be fitted with $\frac{\lambda}{2}$ spacing between any two adjacent antennas. In the case with $W=5$, the benchmark becomes the best antenna selection system with $11$ uncorrelated antennas. In the simulations, our aim is to study how the algorithms perform when the number of observed ports, i.e., $|{\cal K}|$, changes. The parameters used for the LSTM neural network model are listed in TABLE \ref{table:par}.}


\begin{table}[]
\centering
\caption{Parameters for LSTM} 
\begin{tabular}{c|l} 
\thickhline 
Layer & Details \& Values \\ 
\hline\hline 
\multirow{2}{*}{Layer $1$ (input)} & $10$ LSTM cells; {\color{blue}input dimension $(n,1)$}; \\\cline{2-2}
& linear activation\\\hline 
Layer $2$ (hidden) & Dense: $200$ neurons; linear activation \\\hline
Layer $3$ (hidden) & Dropout with probability of $0.2$ \\\hline
Layer $4$ (hidden) & Dense: $200$ neurons; linear activation \\\hline
Layer $5$ (hidden) & Dropout with probability of $0.5$ \\\hline
Layer $6$ (hidden) & Dense: $200$ neurons; linear activation \\\hline
Layer $7$ (hidden) & Dropout with probability of $0.2$ \\\hline
Layer $8$ (hidden) & Dense: $200$ neurons; linear activation \\\hline
Layer $9$ (hidden) & Dropout with probability of $0.5$ \\\hline
Layer $10$ (output) & $N$ neurons; linear activation \\\hline
Batch size & $10$ \\\hline
Number of epochs & $50$ \\\hline
Optimizer & SGD with MSE as loss function \\ [0.1ex] 
\thickhline 
\end{tabular}
\label{table:par}
\end{table}

%

\begin{figure}[!h]
{\color{blue}
\begin{center}
\includegraphics[width=\linewidth]{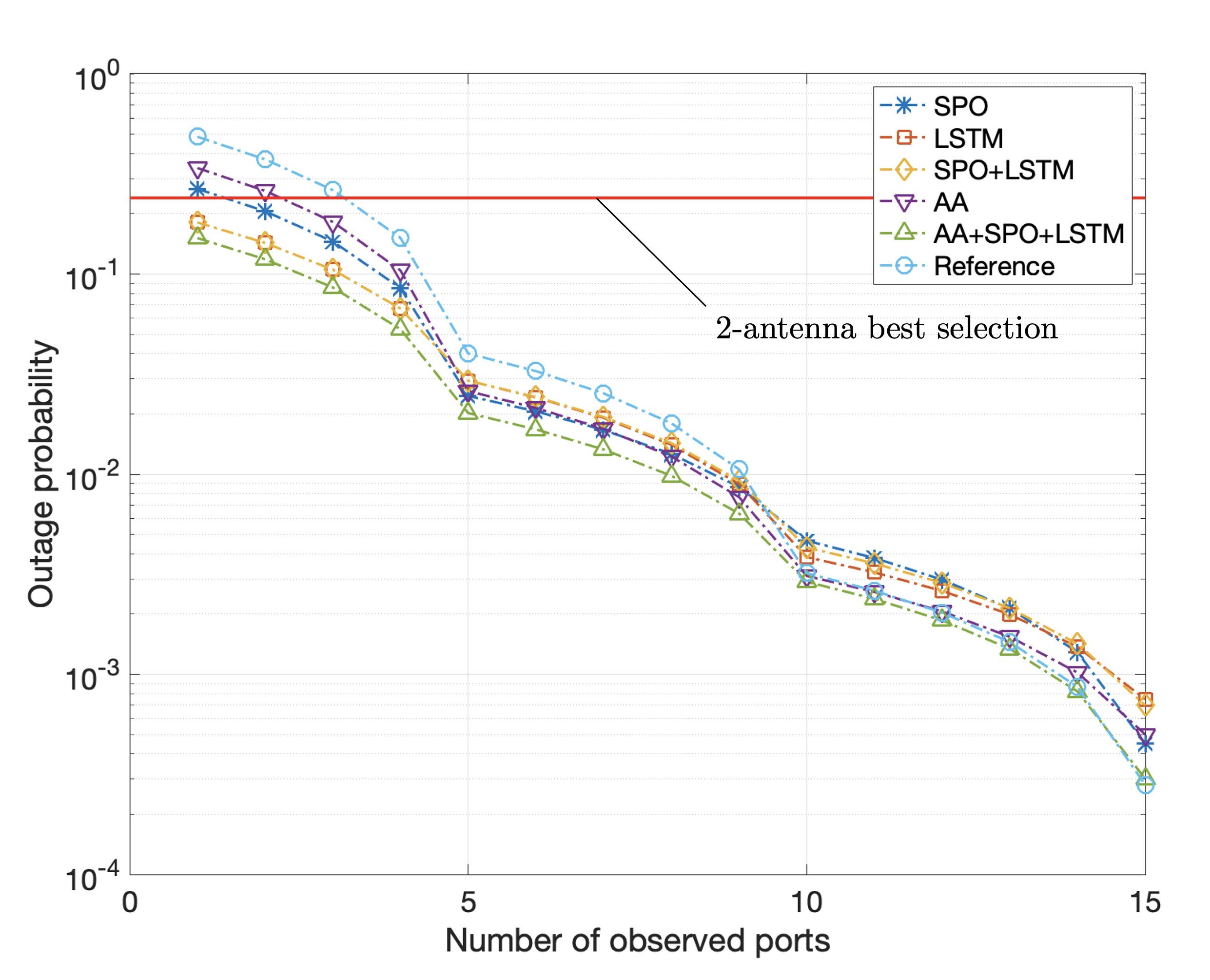}
\caption{Outage probability results when $W=0.5$, $N=50$ and $\gamma=10{\rm dB}$.}\label{fig:Wh}
\end{center}}
\end{figure}

\begin{figure}[!h]
{\color{blue}
\begin{center}
\includegraphics[width=\linewidth]{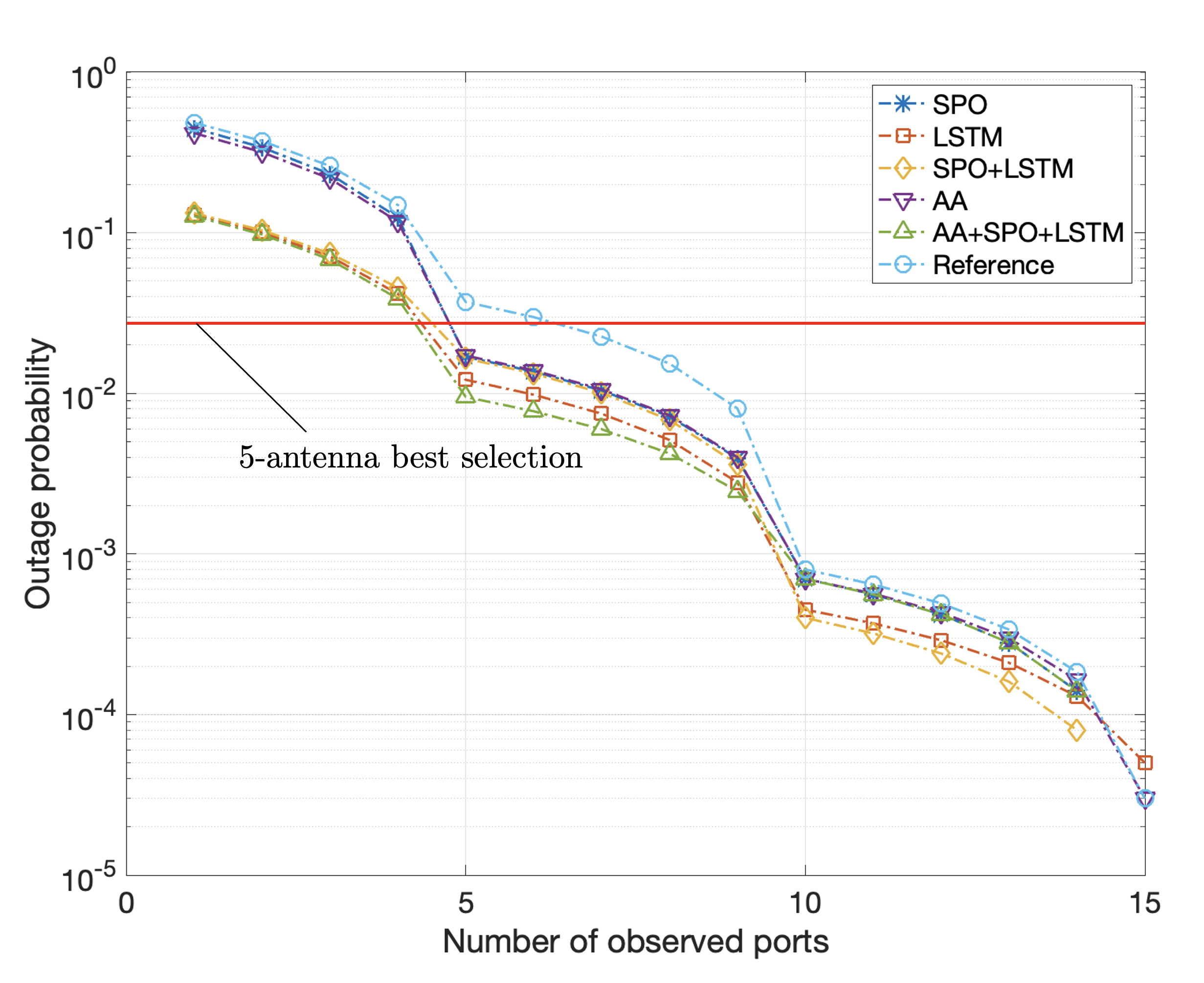}
\caption{Outage probability results when $W=2$, $N=50$ and $\gamma=10{\rm dB}$.}\label{fig:W2}
\end{center}}
\end{figure}

\begin{figure}[!h]
{\color{blue}
\begin{center}
\includegraphics[width=\linewidth]{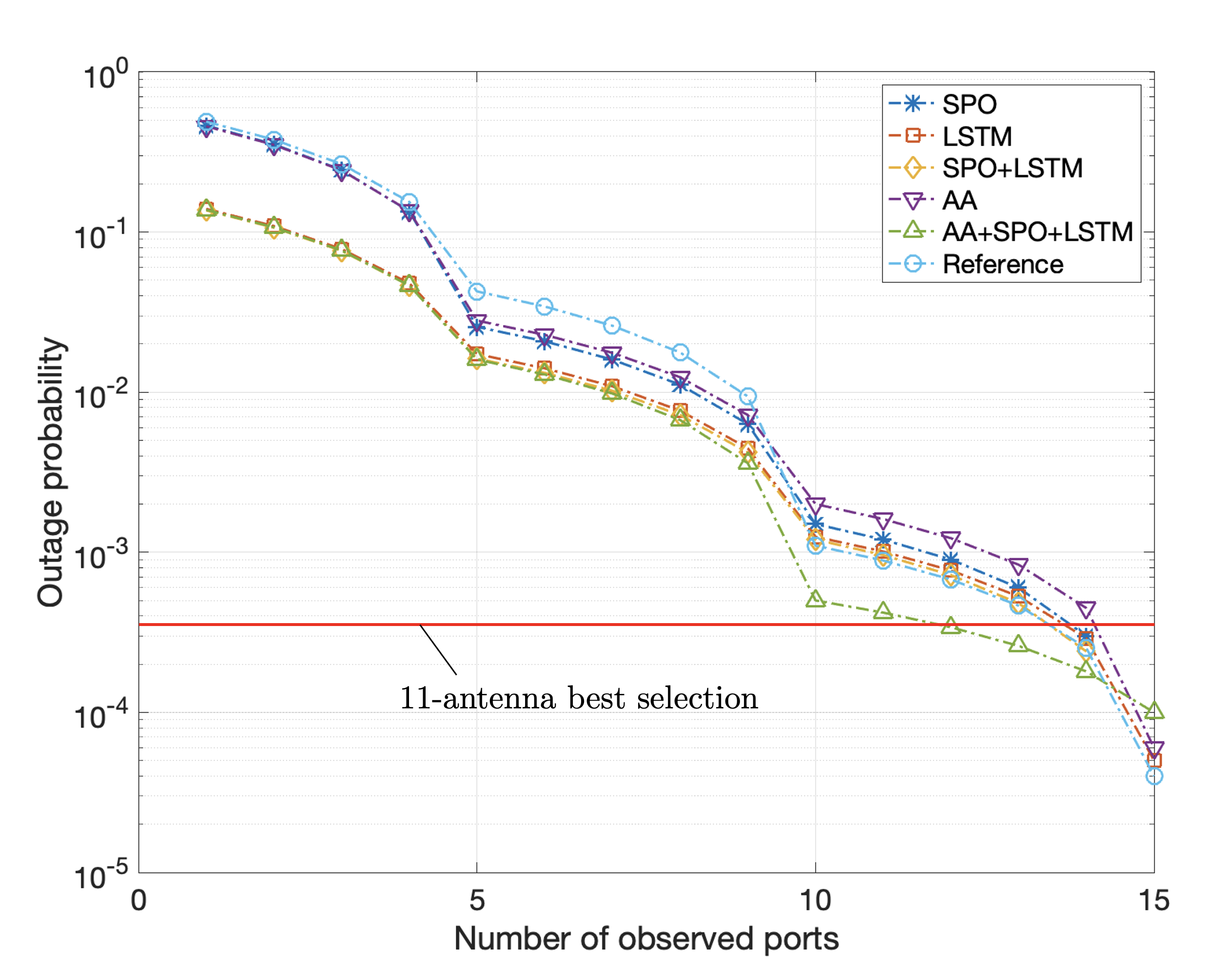}
\caption{Outage probability results when $W=5$, $N=50$ and $\gamma=10{\rm dB}$.}\label{fig:W5}
\end{center}}
\end{figure}

{\color{blue}In Figs.~\ref{fig:Wh}--\ref{fig:W5}, results are provided for outage probability with a target SNR of $10{\rm dB}$ for different values of $W$. The results illustrate that AA+SPO+LSTM performs the best almost in all cases. For $W=2$, if $|{\cal K}|$ gets larger, SPO+LSTM becomes the best but the performance difference from other methods is less apparent. In fact, if $|{\cal K}|$ is large, the performance difference between the methods tends to be smaller. Now, if we compare the results in Figs.~\ref{fig:W2} and $\ref{fig:W5}$, it can be observed that the results for $W=5$ are slightly inferior than that for $W=2$. This is because with a larger space, the observed ports will be farther apart, making the methods more difficult to infer the best port. Additionally, the results show that fluid antenna system with learning-based port selection can outperform significantly the best antenna selection system, with $W=0.5$ indicating the most dramatic gain. This suggests that fluid antenna with practical port selection be a more effective solution to exploit spatial diversity, especially when the size $W$ is small. Finally, the results show that even with only one port observation, the fluid antenna system using AA+SPO+LSTM is able to bring down the outage probability from $0.5$ to close to $0.1$ and if $|{\cal K}|=5$ (i.e., $10\%$ ports observed), more than two orders of magnitude reduction in the outage probability are achieved.}



\section{Conclusion}\label{sec_conclude}
This letter investigated the port selection problem for fluid antenna systems in which only a subset of ports are observed. The problem is crucial to make fluid antenna systems viable. A number of schemes were proposed and {\color{blue}our results showed that significant improvement in outage probability is possible even if very few ports of the fluid antenna system are observed and AA+SPO+LSTM is particularly effective. Our results also indicated that fluid antenna with the proposed port selection algorithms outperform considerably the best antenna selection system with many uncorrelated antennas, which suggests that fluid antenna is able to utilize the space more effectively.}

\ifCLASSOPTIONcaptionsoff
  \newpage
\fi




\begin{thebibliography}{1}

\bibitem{Tse-2003}
L. Zheng, and D. N. C. Tse, ``Diversity and multiplexing: A fundamental tradeoff in multiple-antenna channels,'' {\em IEEE Trans. Inform. Theory}, vol. 49, no. 5, pp. 1073--1096, May 2003.

\bibitem{Larsson-2019}
P. von Butovitsch {\em et al.}, ``Advanced antenna systems for 5G networks,'' [Online] \url{https://www.ericsson.com/en/white-papers/advanced-antenna-systems-for-5g-networks}, White paper, 2019.
\bibitem{Huawei-2019}
``Huawei launches 5G simplified solution,'' [Online] \url{https://www.huawei.com/en/press-events/news/2019/2/huawei-5g-simplified-solution}, 2019.

\bibitem{Wong-2020ell}
K. K. Wong, K. F. Tong, Y. Zhang, and Z. Zheng, ``Fluid antenna system for 6G: When Bruce Lee inspires wireless communications,'' {\em IET Electronics Lett.}, vol. 56, no. 24, pp. 1288--1290, 26 Nov. 2020.
\bibitem{Mitsubishi}
``Mitsubishi electric's SeaAerial antenna uses seawater plume,'' Available [online]: \url{https://www.mitsubishielectric.com}.
\bibitem{Xing-19}
L. Xing, J. Zhu, Q. Xu, D. Yan and Y. Zhao, ``A circular beam-steering antenna with parasitic water reflectors,'' {\em IEEE Antennas and Wireless Propag. Lett.}, vol. 18, no. 10, pp. 2140--2144, Oct. 2019.

\bibitem{Ohta-2013}
A. M. Morishita, C. K. Y. Kitamura, A. T. Ohta, and W. A. Shiroma, ``A liquid-metal monopole array with tunable frequency, gain, and beam steering,'' {\em IEEE Antennas Wireless Propag. Lett.}, vol. 12, pp. 1388--1391, 2013.
\bibitem{Dey-2016}
A. Dey, R. Guldiken, and G. Mumcu, ``Microfluidically reconfigured wideband frequency-tunable liquid-metal monopole antenna,'' {\em IEEE Trans. Antennas Propag.}, vol. 64, no. 6, pp. 2572--2576, Jun. 2016.
\bibitem{Tong-2017}
C. Borda-Fortuny, K.-F. Tong, A. Al-Armaghany, and K. K. Wong, ``A low-cost fluid switch for frequency-reconfigurable Vivaldi antenna,'' {\em IEEE Antennas Wireless Propag. Lett.}, vol. 16, pp. 3151--3154, 2017.
\bibitem{Singh-2019}
A. Singh, I. Goode, and C. E. Saavedra, ``A multistate frequency reconfigurable monopole antenna using fluidic channels,'' {\em IEEE Antennas Wireless Propag. Lett.}, vol. 18, no. 5, pp. 856--860, May 2019.
\bibitem{Shen-2020apmc}
Y. Shen, K. F. Tong, and K. K. Wong, ``Beam-steering surface wave fluid antennas for MIMO applications,'' in {\em Proc. The 2020 Asia-Pacific Microwave Conf. (APMC 2020)}, 8-11 Dec. 2020, Hong Kong, China.

\bibitem{Huang-2021access}
Y. Huang, L. Xing, C. Song, S. Wang and F. Elhouni, ``Liquid antennas: Past, present and future,'' {\em IEEE Open J. Antennas and Propag.}, vol. 2, pp. 473--487, 2021.

\bibitem{Wong-2020twc}
K. K. Wong, A. Shojaeifard, K.-F. Tong and Y. Zhang, ``Fluid antenna systems,'' {\em IEEE Trans. Wireless Commun.}, vol. 20, no. 3, pp. 1950--1962, Mar. 2021.
\bibitem{Wong-2020cl}
K. K. Wong, A. Shojaeifard, K.-F. Tong and Y. Zhang, ``Performance limits of fluid antenna systems,'' {\em IEEE Commun. Lett.}, vol. 24, no. 11, pp. 2469--2472, Nov. 2020.

\bibitem{Convery-2019}
N. Convery, and N. Gadegaard, ``30 years of microfluidics,'' {\em Micro and Nano Engineering}, vol. 2, pp. 76--91, Mar. 2019.


{\color{blue}
\bibitem{Joung-2016}
J. Joung, ``Machine learning-based antenna selection in wireless communications,'' {\em IEEE Commun. Lett.}, vol. 20, no. 11, pp. 2241--2244, Nov. 2016.
\bibitem{He-2018}
D. He, C. Liu, T. Q. S. Quek and H. Wang, ``Transmit antenna selection in MIMO wiretap channels: A machine learning approach,'' {\em IEEE Wireless Commun. Lett.}, vol. 7, no. 4, pp. 634--637, Aug. 2018.}

\bibitem{Grigas-2021}
A. N. Elmachtoub and P. Grigas, ``Smart `predict, then optimize','' {\em Management Science}, 2021, doi: 10.1287/mnsc.2020.3922.

\bibitem{Stuber-2002}
G. L. St$\ddot{\rm u}$ber, {\em Principles of Mobile Communication}, Second Edition, Kluwer Academic Publishers, 2002.

\bibitem{Tao-2021}
Z. Zhang and M. Tao, ``Deep learning for wireless coded caching with unknown and time-variant content popularity,'' {\em IEEE Trans. Wireless Commun.}, vol. 20, no. 2, pp. 1152--1163, Feb. 2021.

\bibitem{Jaggi-2011}
M. Jaggi, ``Convex optimization without projection steps,'' [Online] arXiv preprint arXiv:1108.1170, 2011.

\end{thebibliography}
\end{document}